\documentclass[11pt]{article}
\usepackage{amssymb,amsfonts,amsmath,latexsym,epsf,tikz,url}

\newtheorem{theorem}{Theorem}[section]
\newtheorem{proposition}[theorem]{Proposition}
\newtheorem{conjecture}[theorem]{Conjecture}
\newtheorem{corollary}[theorem]{Corollary}
\newtheorem{lemma}[theorem]{Lemma}

\newtheorem{example}[theorem]{Example}
\newtheorem{definition}[theorem]{Definition}

\newtheorem{problem}[theorem]{Problem}
\newcommand{\proof}{\noindent{\bf Proof.\ }}
\newcommand{\qed}{\hfill $\square$\medskip}

\begin{document}

\title{Distinguishing critical graphs}

\author{
Saeid Alikhani  $^{}$\footnote{Corresponding author}
\and
Samaneh Soltani
}

\date{\today}

\maketitle

\begin{center}
Department of Mathematics, Yazd University, 89195-741, Yazd, Iran\\
{\tt alikhani@yazd.ac.ir, s.soltani1979@gmail.com}
\end{center}


\begin{abstract}
The distinguishing number  $D(G)$   of a graph $G$ is the least integer $d$
such that $G$ has a vertex labeling   with $d$ labels  that is preserved only by a trivial
automorphism. We say that  a graph $G$ is $d$-distinguishing critical, if $D(G)=d$ and $D(H)\neq D(G)$, for every proper induced subgraph $H$ of $G$. 
 This  generalizes the usual definition of a $d$-chromatic critical graph. While the investigation of $d$-critical graphs is a well established
 part of coloring theory, not much is known about $d$-distinguishing  critical graphs. 
In this paper we determine all $d$-distinguishing critical graphs for $d= 1,2,3$ and observe that all of these kind of graphs are $k$-regular graph for some $k\leq d$. 
Also, we show that the disconnected $d$-distinguishing critical graph with $c$ connected
components such that   $c\geq \frac{d}{2}$, is a regular graph.
\end{abstract}


\medskip
\noindent{\bf AMS Subj.\ Class.:} 05C15, 05E18

\section{Introduction}

Let $G=(V,E)$ be a simple graph.  We use the standard graph notation.  In particular, ${\rm Aut}(G)$ denotes the automorphism group of $G$. For simple connected graph $G$, and $v \in V$, the \textit{neighborhood} of a vertex $v$ is the set $N_G(v) = \{u \in V(G) : uv \in   E(G)\}$. The \textit{degree} of a vertex $v$ in a graph $G$, denoted by ${\rm deg}_G(v)$, is the number of edges of $G$ incident with $v$. In particular, ${\rm deg}_G(v)$ is the number of neighbours of $v$ in $G$.  We denote by $\delta(G)$ and $\Delta(G)$ the minimum and maximum degrees of the vertices of $G$. A graph $G$ is \textit{$k$-regular} if ${\rm deg}_G(v) = k$ for all $v \in V$. 
For $u, v \in V$ the \textit{distance} of
 $u$ and $v$, denoted by $d_G(u, v)$ or $d(u, v)$, is the length of a shortest path between
 $u$ and $v$. The \textit{eccentricity}  of a vertex $v$, ${\rm ecc}(v)$, is ${\rm max}\{d(u, v) : u \in V \}$. The maximum and minimum eccentricity of vertices of $G$ are called \textit{diameter}  and
 \textit{radius} of $G$ and are denoted by ${\rm diam}(G)$ and ${\rm rad}(G)$, respectively. \textit{Center}  of $G$ is the subgraph induced by vertices with eccentricity ${\rm rad}(G)$. A graph is
 called \textit{self-centered}  if it is equal to its center, or equivalently, its diameter equals
 its radius.  A graph $G$ is called  \textit{$k$-self-centered} if ${\rm diam}(G) = {\rm rad}(G) = k$. The terminology \textit{k-equi-eccentric}  graph is also used by some authors. For studies on these
 graphs see \cite{Buckly}.
A \textit{clique}  of a graph is a set of mutually adjacent vertices, and that the maximum
size of a clique of  $G$, the clique number of $G$, is denoted $\omega(G)$. 
The complementary notion of a clique is an independent set, a set of vertices no two of
which are adjacent. An independent set in a graph is maximum if the graph contains no
larger stable set. The cardinality of a maximum independent set in  $G$ is called the
 \textit{independence number} of $G$, denoted by $\alpha (G)$. Clearly, a subset $S$ of $V$ is an independent set in $G$ if and only if $S$ is a clique in $\overline{G}$, 
 the complement of $G$. In particular, $\omega(\overline{G}) = \alpha(G)$.

A labeling of $G$, $\phi : V \rightarrow \{1, 2, \ldots , r\}$, is said to be \textit{$r$-distinguishing}, 
if no non-trivial  automorphism of $G$ preserves all of the vertex labels.
The point of the labels on the vertices is to destroy the symmetries of the
graph, that is, to make the automorphism group of the labeled graph trivial.
Formally, $\phi$ is $r$-distinguishing if for every non-trivial $\sigma \in {\rm Aut}(G)$, there
exists $x$ in $V = V (G)$ such that $\phi(x) \neq \phi(\sigma(x))$. We will often refer to a
labeling as a coloring, but there is no assumption that adjacent vertices get
different colors. Of course the goal is to minimize the number of colors used.
Consequently  the \textit{distinguishing number} of  $G$ is defined  by

\begin{equation*}
D(G) = {\rm min}\{r \vert ~ G ~\textsl{has a labeling that is $r$-distinguishing}\}.
\end{equation*} 

This number has introduced in \cite{Albert}.  If a graph has no nontrivial automorphisms, its distinguishing number is  $1$. In other words, $D(G) = 1$ for the asymmetric graphs.
 The other extreme, $D(G) = \vert V(G) \vert$, occurs if and only if $G = K_n$. The distinguishing number of some examples of graphs was obtained in \cite{Albert}. For 
 instance, $D(P_n) =2$ for every $n\geq 3$, and  $D(C_n) =3$ for $n =3,4,5$,  $D(C_n) =2$ for $n \geq 6$. Also, $D(K_{n,m})= n$ where $n > m\geq 1$, and $D(K_{n,n}) =n+1$ for any $n \geq 3$.  Albertson and Collins in \cite{Albert} showed that $D(G)=D(\overline{G})$, where $\overline{G}$ is the complement  graph of $G$.
 
 \medskip

 A graph $G$ is (colour) critical if $\chi(H)<\chi(G)$ for every proper subgraph $H$ of $G$, and $G$ is said to be $k$-critical if $G$ is critical and $\chi(G)=k$.
  The investigation of $k$-critical graphs is a well established   part of coloring theory. Actually, criticality is a general concept in graph theory and can be defined with respect to various graph parameters or graph  properties. The importance of the notion of criticality lies in the fact that problems for graphs may often be reduced to
 problems for critical graphs, and the structure of the latter is more restricted. Critical graphs with respect to the chromatic
 number were first defined and used by Dirac \cite{dirac}  in 1951  (\cite{list}).
 Stiebitz, Tuza and Voigt introduced and discussed the list critical graphs in \cite{list}. 
 
 \medskip 
 
 In this paper we introduce  critical graphs with respect to the distinguishing number and  we discuss some basic properties of $d$-distinguishing critical graphs. 
 More precisely, we characterize connected $d$-distinguishing critical graphs for $d=1,2$   in Section 2. We state some necessary condition for $3$-distinguishing critical graphs in this section, too. In Section 3, we study  disconnected $d$-distinguishing critical graphs and show that there are exactly five $3$-distinguishing critical graphs. Also we obtain disconnected  $d$-distinguishing graphs for $d=5,6$. Finally we propose two conjectures and a problem in the last section.   

\section{Characterization of  $d$-distinguishing critical graphs }

Criticality is a general concept in graph theory and can be defined with respect to various graph parameters or graph  properties. Motivated by chromatic critical graphs and 
list critical graphs we state the following definition.

\begin{definition}
A graph $G$ is {\rm $d$-distinguishing critical}, if $D(G)=d$ and $D(H)\neq D(G)$, for every proper induced subgraph $H$ of $G$.
\end{definition}

Let us first discuss some elementary facts about distinguishing  critical graphs by the following example. 

\begin{example}~
\begin{itemize}
\item The complete graph  $K_n$ (and also its complement) is an $n$-distinguishing critical graph for every $n \geq 3$. 
\item The complete bipartite graph  $K_{n,n}$ (and also its complement) is a $(n+1)$-distinguishing critical graph for every $n \geq 3$. 
\item The cycles graphs $C_3, C_4$ and $C_5$ are 3-distinguishing critical graphs.
\end{itemize}
\end{example}

Since the distinguishing number of each graph and its complement is equal, so we have the following proposition:  

\begin{proposition}\label{compdiscrt}
A graph $G$ is $d$-distinguishing critical, if and only if $\overline{G}$ is a $d$-distinguishing critical graph.
\end{proposition}

\begin{figure}[h]
	\begin{center}
		\includegraphics[width=0.85\textwidth]{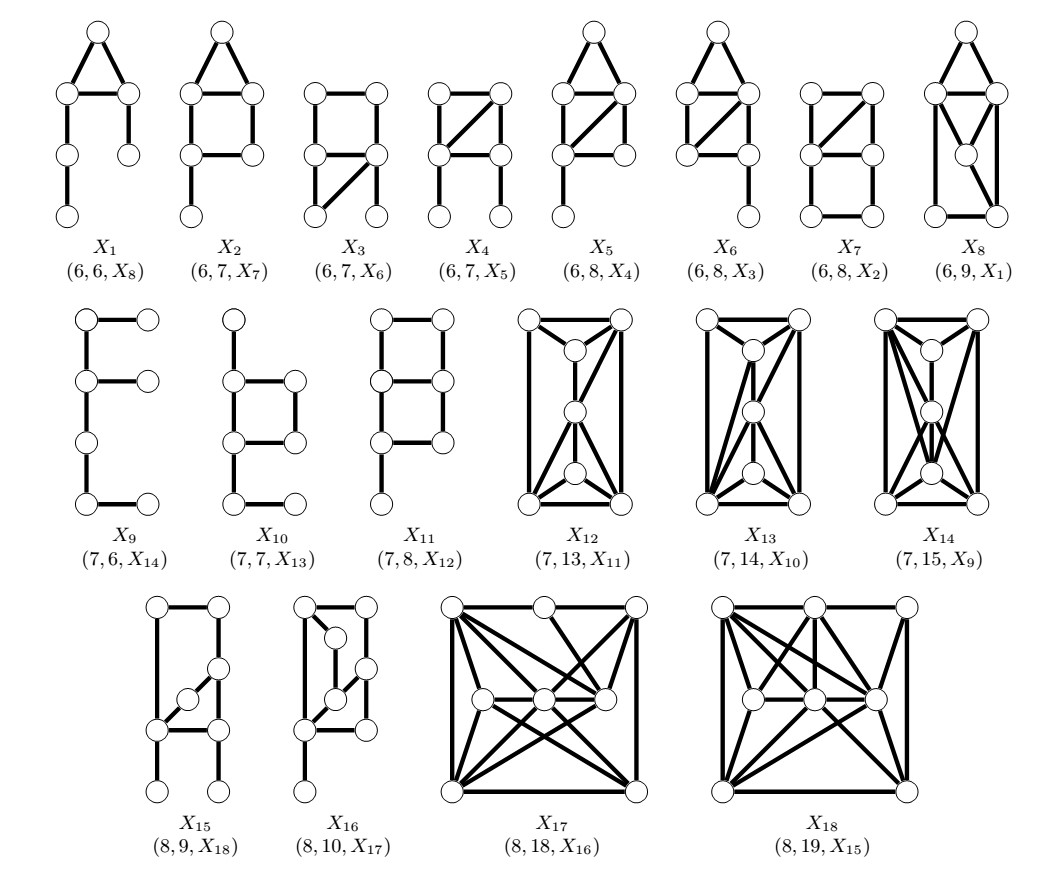}
		\caption{ \label{fig1new} \small The $18$ minimal asymmetric graphs. 
For each graph the triple $(n,m, co-G)$, describes the number of vertices, edges and the name
of the complement graph, respectively. The graphs are ordered first by number of vertices and
second by number of edges.}
	\end{center}
\end{figure}

First, we determine  $d$-distinguishing critical graphs for all $d\leq 2$. An undirected graph $G$ on at least
two vertices is \textit{minimal asymmetric} if $G$ is asymmetric and no proper induced subgraph of $G$ on at least two vertices is asymmetric.  Pascal Schweitzer and Patrick Schweitzer in \cite{Schweitzer} showed that there are exactly $18$ finite minimal asymmetric undirected graphs up to isomorphism. These $18$ graphs depicted in Figure \ref{fig1new}.
\begin{theorem}~
	\begin{enumerate} 
		\item[(i)] There is no  $1$-distinguishing critical graph of order $n\geq 2$.
		
		\item[(ii)] The only $2$-distinguishing critical graphs are $K_2$ and $\overline{K_2}$.
		\end{enumerate} 
		\end{theorem}
\proof 
	\begin{enumerate} 
		\item[(i)]
		Since $K_1$ is a proper induced subgraph of $G$ with the distinguishing number $1$, so the result follows.
		
\item[(ii)] Let $G$ be a 2-distinguishing critical graph of order $n\geq 3$. Then there exist at least two distinct vertices of $G$, say $v$ and $w$. If $H$ is the proper induced subgraph of $G$ generated by $v$ and $w$, then $D(H)=2$, and hence $G$ is not a $2$-distinguishing critical graph, which is a contradiction. \qed
\end{enumerate}

In sequel, we want to characterize the  $d$-distinguishing critical graphs with $d \geq 3$. We need the following theorem. 

 \begin{theorem}{\rm \cite{EJC,wongetall}}\label{uppboundisnumb}
 	If $G$ is a connected graph with maximum degree $\Delta$, then $D(G) \leq \Delta +1$. Furthermore, equality holds if and only if $G$ is a $K_n$, $K_{n,n}$ or $C_5$.
 \end{theorem}

\begin{theorem}\label{characdiscritgra}
Let  $G$ be a  $d$-distinguishing critical graph with $d \geq 3$ and order $n$. Then the following properties are satisfied:
\begin{itemize}
\item[{\rm (i)}] If   $G$ is connected and $\Delta (G) = d-1$, then $G$ is $K_d, K_{d-1,d-1}$ or $C_5$.
\item[{\rm (ii)}] If  $G$ is connected  and $G\neq K_d , K_{d-1,d-1}, C_5$, then $\Delta (G) \geq d$.
\item[{\rm (iii)}] If  $G$ and $\overline{G}$ are connected and $G\neq C_5$, then $\delta (G) \leq n-d-1$.
\item[{\rm (iv)}] The graph $G$ is a $K_{1,d}$-free graph, i.e., $G$ does not have  the star graph $K_{1,d}$ as  its induced subgraph. 
\item[{\rm (v)}]  If $G\neq K_d$, then the clique number of $G$ is  $\omega(G) \leq d-1$.
\item[{\rm (vi)}] If $G \neq \overline{K_d}$, then the independence number of $G$ is $\alpha (G) \leq d-1$.
\end{itemize}
\end{theorem}
\proof  The cases  (i) and (ii) follow directly from Theorem \ref{uppboundisnumb}. It can be concluded from case (ii) that $n-1- \delta (G) = \Delta (\overline{G}) \geq d$, and this prove the case (iii).  
To prove the case (iv), it is sufficient to know that $D(K_{1,d})=d$. For (v), if $\omega(G) \geq d$, then $G$ has a proper induced subgraph $K_d$ with $D(K_d)=d$, which is a contradiction. Finally, if $\alpha (G) \geq d$, then the complement graph $\overline{K_d}$ is a proper induced subgraph of $G$ with the distinguishing number $d$, which is a contradiction. \qed

\begin{theorem}\label{triclawcrit}~
\begin{itemize}
 \item[\rm{(i)}]   Let $G$ be a connected $d$-distinguishing critical graph for some $d\geq 2$. If  $G$ is a triangle free graph, then $G= K_2, K_{d-1,d-1}$, or $C_5$.
\item[\rm{(ii)}] Let $G$ and $\overline{G}$ be two connected graphs. If $G$ is a claw free graph such that $G\neq C_5$, then $G$ is not a $d$-distinguishing critical graph for any $d\geq 2$.
\end{itemize}
\end{theorem}
\proof  
\begin{enumerate}
 \item[(i)]
 Let $G$ be a connected $d$-distinguishing critical graph for some $d\geq 2$, such that $G\neq K_2, K_{d-1,d-1}, C_5$. We show that $G$ has a triangle. Let $v$ be a vertex of $G$ with maximum degree $\Delta$, and $N_G(v)=\{v_1, \ldots , v_{\Delta}\}$. By Theorem \ref{characdiscritgra}, we have $\Delta \geq d$ and $\alpha(G) \leq d-1$, so there exist at least two adjacent vertices in $N_G(v)$, say $v_1$ and $v_2$. Hence, the graph generated by vertices $v, v_1, v_2$ makes a triangle graph as the induced subgraph of $G$.

\item[(ii)] By contradiction suppose that  $G$ is a $d$-distinguishing  critical graph for some $d\geq 2$. So $\overline{G}$ is a $d$-distinguishing  critical graph, by Proposition \ref{compdiscrt}. Since $G$ is claw free, so $\overline{G}$ is a triangle free graph. Now by Part (i), and connectivity of  $G$ and $\overline{G}$ we conclude that $G= C_5$, which is a contradiction. \qed
\end{enumerate}

\section{Disconnected distinguishing critical graphs}

A vertex $v$ in a graph $G$ of order $n$, is called an $a_1$-vertex if ${\rm deg}_G(v) = n - 1$
or there is a vertex $u \in N_G(v)$ such that $N_G(u)\cup  N_G(v) = V (G)$; otherwise, $u$ is an
$a_2$-vertex, see \cite{Malaravan}. Malaravan in \cite{Malaravan}, showed that if ${\rm diam}(G) \geq 4$, then every vertex of $G$ is an $a_2$-vertex. 
\begin{theorem}{\rm \cite{Malaravan}}\label{selfcenter}
Let  $G$ be  a graph with no isolates. Then, each vertex of $G$ is an
$a_2$-vertex in $G$ if and only if $\overline{G}$ is a 2-self-centered graph. 
\end{theorem}
\begin{theorem}
 If $G$ is a disconnected  $d$-distinguishing critical  graph with no isolated vertex, then $\overline{G}$ is a connected 2-self-centered graph.
\end{theorem}
\proof  Since $G$ is a disconnected graph, so ${\rm diam}(G) \geq 4$, and hence every vertex of $G$ is an $a_2$-vertex. Now, we have the result by Theorem \ref{selfcenter}.\qed

Let $(G, \phi)$ denote the labeled version of $G$ under the labeling $\phi$. Given two distinguishing $k$-labelings $\phi$ and $\phi'$ of $G$, we say that $\phi$ and $\phi'$ are \textit{equivalent} if there is some automorphism of $G$ that maps $(G, \phi)$ to $(G, \phi')$. We need $D(G, k)$, the number of inequivalent $k$-distinguishing labelings of $G$, which was first considered by Arvind and Devanur \cite{Arvind and Devanur} and Cheng \cite{ChengD(Gk)} to determine $d$-distinguishing critical graph with $d \geq 3$.

\begin{lemma}{\rm \cite{Arvind and Devanur}}\label{rem1}
Let $H$ be a connected graph. If $G$ consists of $c(G)$ copies of $H$, i.e., $G= c(G) H$, then
$D(G) = {\rm min}\{k : D(H, k) \geq c(G)\}$.
\end{lemma}

\begin{lemma}\label{lemineqcritic}
 Let $G$ be a disconnected $d$-distinguishing critical graph ($d\geq 3$) with $c(G)$ connected  components.
\begin{itemize}
\item[{\rm (i)}] The connected components of $G$ are pairwise isomorphic.
\item[{\rm (ii)}] If the connected components of $G$ are isomorphic to graph $H$, then
\begin{enumerate}
\item[1)]  $D(n H) < d$ for every $1\leq n < c(G)$, and 
\item[2)]  $c(G)= D(H, d-1) +1$, and $D(H, d) \geq c(G)$, and
\item[3)] $c(G) > D(H, D(H))$.
\end{enumerate}
\end{itemize}
\end{lemma}
\proof 
\textbf{(i)} Let $G_{i1}, \ldots , G_{is_i}$ be the connected components of $G$ such that $G_{i1} \cong \cdots \cong  G_{is_i}$ for every  $1\leq i \leq k$, and also $G_{i1}\ncong G_{j1}$ for $i\neq j$.  Set $H_i = G_{i1}\cup \cdots \cup G_{is_i}$.  It is clear that  $D(G) = {\rm max} \{D(H_i)\}_{i=1}^k$. Without loss of generality we can assume that $D(H_1) = {\rm max} \{D(H_i)\}_{i=1}^k$. If $k\geq 2$, then the distinguishing number of  the proper induced subgraph $H_1$ of $G$ is $D(G)$, which is a contradiction. Therefore $k=1$, i.e., all connected components of $G$ are isomorphic to each other.

\medskip
\textbf{(ii-1)}  By Lemma \ref{rem1}, we have $D(nH) \leq d$ for every $1\leq n < c(G)$. If there exist $n$ such that  $D(n H) = d$, then thr graph $K:= n H$ is a proper induced subgraph of $G$ with $D(K) = D(G)$, which is a contradiction.

\textbf{(ii-2)} By part (1), and Lemma \ref{rem1},  we have $D((c(G)-1)H)\leq d-1$ and $D(H, d-1) < c(G)$. Thus $D(H, d-1)\geq c(G)-1$ and $D(H, d-1) < c(G)$. Therefore, $c(G)= D(H, d-1) +1$.

\textbf{(ii-3)} By  contradiction suppose that  $c(G) \leq D(H, D(H))$. So $d= D(G)= D(H)$, by Lemma \ref{rem1}. Since $c(G)\geq 2$, so $H$ is a proper induced subgraph of $G$ with $D(H)= D(G)$, which is a contradiction.\qed

Before we prove the next result we need some additional information  about the distinguishing number of  complete multipartite graphs. Let $K_{{a_1}^{j_1},{a_2}^{j_2},\ldots, {a_r}^{j_r}}$  denote the complete multipartite graph that has $j_i$ partite sets of size $a_i$ for $i = 1, 2,\ldots , r$ and $a_1 > a_2 > \ldots > a_r$.

\begin{theorem}\label{Theorem 2.4 of Collins and A. N. Trenk}{\rm \cite{EJC}}
	Let $K_{{a_1}^{j_1},{a_2}^{j_2},\ldots, {a_r}^{j_r}}$ denote the complete multipartite graph that has $j_i$ partite sets of size $a_i$ for $i = 1, 2,\ldots , r$, and $a_1 > a_2 > \ldots > a_r$. Then 
	\begin{equation*}
	D(K_{{a_1}^{j_1},{a_2}^{j_2},\ldots, {a_r}^{j_r}}) ={\rm min}\left\{p : {p \choose a_i} \geq  j_i ~{\rm for~all~i}\right\}.
	\end{equation*}
\end{theorem}

Since the distinguishing number of a graph and its complement is equal, so if $G$ consists of $c(G)$ copies of complete graph $K_s$, then $D(c(G) K_s)= {\rm min}\{p : {p \choose s} \geq c(G) \}$.

\begin{theorem}\label{caracterdlessd2}
 Let $G$ be a disconnected $d$-distinguishing critical graph ($d\geq 3$) with $c(G)$ isomorphic connected  components $H$, and the independence number $\alpha (G)$.
\begin{enumerate}
\item[{\rm (i)}] If $c(G)\geq \frac{d}{2}$, then $G = c(G) K_s$, where  $c(G)= {d-1 \choose s}+1$ and ${d \choose s}\geq c(G)$.
\item[{\rm (ii)}]  If $\alpha (G) = c(G)$, or $\alpha (G)$ is a prime number, then  $G = c(G) K_s$, where  $c(G) ={d-1 \choose s} +1$ and ${d \choose s}\geq c(G)$.
\end{enumerate}
\end{theorem}
\proof
\begin{enumerate} 
\item[(i)]  We consider the two following cases:

Case 1)  If $c(G) = d$, then by choosing one vertex from each component of $G$ we can make the induced subgraph  $\overline{K_d}$ of $G$. Since $G$ is a $d$-distinguishing critical graph, so $G= \overline{K_d}$.

Case 2) If $c(G) \leq d-1$, then $G\neq \overline{K_d}$, and so $\alpha (G) \leq d-1$, by part (vi) of Theorem \ref{uppboundisnumb}. Since $\alpha (G)= c(G). \alpha (H)$ and $d\leq 2c(G)$, so we obtain $\alpha (H) = 1$, and hence $H$ is a complete graph, say $K_s$. Thus $G= c(G)K_s$. Now by $D(G)= d$ and $D(c(G)K_s)= {\rm min}\{p~:~{p\choose s}\geq c(G)\}$ (by Theorem \ref{Theorem 2.4 of Collins and A. N. Trenk}), we have 
\begin{equation}\label{eq1}
d={\rm min}\left\{p~:~{p\choose s}\geq c(G)\right\}.
\end{equation}

On the other hand the distinguishing number of the proper induced subgraph  $H= (c(G)-1) K_s$ is $D(H)\leq d-1$, by Lemma \ref{lemineqcritic} (i). In fact  $D(H)= d-1$, since otherwise $D(H)\leq d-2$, and so $D(H\cup K_s) = D(G) \leq d-1$, which is a contradiction. Thus 
\begin{equation}\label{eq2}
d-1={\rm min}\left\{p~:~{p\choose s}\geq c(G)-1\right\}.
\end{equation}

Now using Equations \eqref{eq1} and \eqref{eq2}, we conclude that ${d-1 \choose s} = c(G)-1$ and ${d \choose s}\geq c(G)$.

\item[(ii)]   If $\alpha (G) = c(G)$, or $\alpha (G)$ is a prime number, then we obtain $\alpha (H) =1$, from   $\alpha (G)= c(G) \alpha (H)$. Hence $H$ is a complete graph, say $K_s$. The remaining proof is the same as the proof of  case (2) of Part (i). \qed
\end{enumerate} 

To determine exactly $3$-distinguishing critical graphs we need the following lemma.
\begin{lemma}\label{lemm}
 If $G$ is a $3$-distinguishing critical graph, then $\Delta (G) \leq 2$. 
\end{lemma}
\proof  By   contradiction, suppose  that $G$ has a vertex $v$ of degree at least three. Let  $w,y,z$ be  three distinct adjacent vertices to $v$. If at least the two vertices of $w,y,z$ are adjacent, say $w,y$, then the proper induced subgraph $H$ of $G$ generated by $v,w,y$ is a triangle with the distinguishing number 3, which is a contradiction. If non of vertices  $w,y,z$ are adjacent,  then the proper induced subgraph $H$ of $G$ generated by $v,w,y,z$ is the star graph $K_{1,3}$  with the distinguishing number $3$, which is a contradiction. \qed 

\begin{theorem}
There are exactly five $3$-distinguishing critical graphs, $C_3, C_4, C_5, \overline{C_3}, \overline{C_4}$. 
\end{theorem}
\proof By Theorem \ref{uppboundisnumb} and Lemma \ref{lemm}   it is easy to see that the  graphs  $C_3, C_4, C_5$ are all connected $3$-distinguishing critical graph. Also  the complements of graphs $C_3, C_4, C_5$ are the only disconnected  $3$-distinguishing critical graphs, by  Theorem \ref{caracterdlessd2} (i).\qed


\begin{theorem}\label{lemdisconclessd}
If $G$ is a  disconnected $d$-distinguishing critical graph ($d\geq 5$) with $c(G)$ isomorphic connected  components $H$,  such that $c(G) < \frac{d}{2}$, then $d= D(H)+1$, $D(H) \geq 4$ and $D(H, D(H)) < \frac{D(H)-1}{2}$.
\end{theorem}
\proof  By part (ii-1) of Lemma \ref{lemineqcritic}, we can suppose that $d = D(H)+i$  for some $i \geq 1$.  By Lemma \ref{rem1}, we have $c = D(H, D(H)+i-1)+1$. Since there exist ${D(H)+i-1 \choose D(H)}$ distinct sets  with $D(H)$ labels,  so ${D(H)+i-1 \choose D(H)} \leq D(H, D(H)+i-1)$. Now since $c(G) < \frac{d}{2}$, so we can conclude that 

\begin{equation}\label{eqi}
 {D(H)+i-1 \choose D(H)} < \frac{D(H)+i-2}{2}.
\end{equation}
It can be seen that if $i\geq 2$, then by Equation \eqref{eqi} we have  $D(H)\leq 0$, which is a contradiction. Thus $i =1$, and so we obtain $D(H) \geq 4$ from  Equation \eqref{eqi}. Also $c(G) = D(H, D(H))+1$. The last inequality  follows directly from $c(G) < \frac{d}{2}$.\qed

\begin{corollary}
If $G$ is a  disconnected $d$-distinguishing critical graph ($d\geq 5$) with $c(G)$ isomorphic connected  components $K_s$ (the complete graph of order $s$),  such that $c(G) < \frac{d}{2}$, then $s=d-1$, and $G= 2K_{d-1}$.
\end{corollary}
\proof  It is sufficient to know that $D(K_s , D(K_s)) = D(K_s , s) =1$. Now the result follows from 
Theorem \ref{lemdisconclessd}.\qed

\begin{lemma}\label{claim3}
Let $G$ be a  disconnected $d$-distinguishing critical graph ($d\geq 5$) with $c(G)$ isomorphic connected  components $H$ such that $H$ is not a complete graph, and also $c(G) < \frac{d}{2}$. Let $\phi$ be an arbitrary distinguishing labeling of $H$ with $D(H)$ labels and  $X_i = \{x \in V(H):~\phi(x) = i\}$ where  $i= 1, 2, \ldots , D(H)$. If the  degrees sequence  of vertices in $X_i$  is denoted by the sequence $S_i$, then for every $i,j \in \{1, 2, \ldots , D(H)\}$ we have $S_i = S_j$.
\end{lemma}
\proof Let $n_i$ be  the number of vertices of degree $p$ in $X_i$. By contradiction  assume  that there exists a degree $p$ for which $B_p\neq \emptyset$ where $A_p = \{X_k: n_k = n_1 \}$ and  $B_p = \{X_k:n_k \neq n_1\}$.  
We consider the two following cases:

\medskip
Case 1) If $|A_p| < \frac{D(H) - 1}{2}$, then $|B_p| \geq  \frac{D(H) - 1}{2}$. For every $X_k \in B_p$, we define the distinguishing labeling $\varphi_k$ of $H$ as follows:
\begin{equation*}
\varphi_k(x)=\left\{
\begin{array}{ll}
\phi (x)& x\notin X_1, X_k,\\
1& x \in X_k,\\
k & x\in X_1.
\end{array}\right.
\end{equation*}
Since the number of vertices of degree $p$ in $X_k$ is not equal with the number of vertices of degree $p$ in $X_1$,  we conclude that the distinguishing labelings $\varphi_k$,
 $X_k$, are nonisomorphic distinguishing labeling to each other and also to $\phi$. Since $|B| \geq  \frac{D(H) - 1}{2}$, so $D(H, D(H)) \geq \frac{D(H)-1}{2}$, which is a contradiction  by Theorem \ref{lemdisconclessd}. 
 
 \medskip
 Case 2) If $|A_p| \geq \frac{D(H) - 1}{2}$.  By a similar reasoning as Case 1, we obtain a contradiction.\qed

\begin{theorem} \label{degreenumber}
Let $G$ be a  disconnected $d$-distinguishing critical graph ($d\geq 5$) with $c(G)$ isomorphic connected  components $H$ such that $H$ is not a complete graph, and also $c(G) < \frac{d}{2}$. If $\{p_1, \ldots , p_s\}$ is the set of degrees of vertices of $G$, then the number of vertices of degree $p_i$ in $G$ is at least $2 c(G) (d-1)$, for any $1\leq i \leq s$.
\end{theorem}
\proof   Let $\phi$ be a  distinguishing labeling of $H$ with $D(H)$ labels and  $X_i = \{x \in V(H):~\phi(x) = i\}$ where  $i= 1, 2, \ldots , D(H)$. It is clear that  $\{p_1, \ldots , p_s\}$ is the set of degrees of vertices of $H$, too.  By Lemma \ref{claim3}, we  know that the number of verices of degree $p_i$ in each of sets $X_1, \ldots , X_{D(H)}$ are equal, for any $1\leq i \leq s$. Let $t_i$ is the number of vertices of degree $p_i$ in each of sets $X_1, \ldots , X_{D(H)}$. 
First, we prove that  $t_i \geq 2$ for  any $1\leq i \leq s$. By contradiction, we assume that there exists $p \in \{p_1, \ldots , p_s\}$ such that the number of vertices of degree $p$ in each of sets $X_1, \ldots , X_{D(H)}$ is 1.  So we have exactly $D(H)$ vertices of degree $p$ in graph $H$, say $v_1, \ldots , v_{D(H)}$ such that $v_i \in X_i$ where $1\leq i \leq D(H)$. We claim that the restriction of the automorphism group of $H$ to the set $\{v_1, \ldots , v_{D(H)}\}$ is isomorphic to the permutation group $\mathbb{S}_{D(H)}$. 

To prove our claim note that for every $a,b \in \{1, \ldots , D(H)\}$, we have a labeling $\phi_{a,b}$ of $V(H)$ with $X^{a,b}_i = \{x \in V(H):~\phi_{a,b}(x) = i\}$ 
such that $X^{a,b}_a = X_a - \{v_a\}$, $X^{a,b}_b = X_b \cup \{v_a\}$, and $X^{a,b}_i = X_i$ for every $i \in \{1, \ldots , D(H)\}\setminus \{a,b\}$.  Since the number of vertices of degree $p$ in $X^{a,b}_a$ and $X^{a,b}_b$ are distinct, so the labeling $\phi_{a,b}$ is not distinguishing, by Lemma \ref{claim3}. Hence, there exists a nonidentity automorphism $f_{a,b}$ of $H$ such that $f_{a,b}$ preserves the labeling  $\phi_{a,b}$. Thus $f_{a,b} (v_a) = v_b$ and $f_{a,b} (v_b) = v_a$, since otherwise  $f_{a,b} (v_b) = v_b$ and  $f_{a,b} (v_a) = v_a$, and so  $f_{a,b}$ preserves the distinguishing labeling $\phi$, and hence  $f_{a,b}$ should be the identity automorphism, which is a contradiction. Therefore $f_{a,b} (v_a) = v_b$,  $f_{a,b} (v_b) = v_a$, and $f_{a,b} (v_i) = v_i$ for every $i \in \{1, \ldots , D(H)\}\setminus \{a,b\}$. Thus for any   $a,b \in \{1, \ldots , D(H)\}$, there exists an automorphism $f_{a,b}$ of $H$ such that the restriction of $f_{a,b}$ to $v_1, \ldots , v_{D(H)}$ is isomorphic to the transposition $(a,b)$ in  $\mathbb{S}_{D(H)}$. Hence, we can conclude that the restriction of automorphism group of $H$ to the set $\{v_1, \ldots , v_{D(H)}\}$ is isomorphic to the permutation group $\mathbb{S}_{D(H)}$, and this proves the claim.

 Now since $H$ is a connected graph, so  by claim the set of vertices $\{v_1, \ldots , v_{D(H)}\}$  generates the induced subgraph $K_{D(H)}$ of $H$.  Since $c(G)\geq 2$, so we can make the induced subgraph $K= K_{D(H)} \cup K_{D(H)}$ of $G$ (here $\cup$ is the disjoint union symbol) with $D(K) = D(H)+1$. On the other hand, $d = D(H) +1$ by Theorem \ref{lemdisconclessd}, and so $D(K) =d$. Since $H$ is not a complete graph, so $K$ is a proper induced subgraph of $G$ with $D(K) = D(G)$, which is contradiction to that $G$ is a $d$-distinguishing critical graph. Therefore $t_i \geq 2$, for any $1\leq i \leq s$. 
 
 To complete the proof, it is sufficient to know that the number of vertices of degree $p_i$ in $G$ is $c(G)D(H)t_i$ for any $1\leq i \leq s$, and also $D(H)= d-1$.\qed

\begin{theorem}\label{regcrit}
If $G$ is a  disconnected $d$-distinguishing critical graph ($d\geq 5$) with $c(G)$ isomorphic connected  components $H$,  such that $c(G) < \frac{d}{2}$, then $D(H) ~|~|V(H)|$.
\end{theorem}
\proof   Let $\phi$ be  an arbitrary distinguishing labeling of $H$ with $D(H)$ labels and  $X_i = \{x \in V(H):~\phi(x) = i\}$ where  $i= 1, 2, \ldots , D(H)$. By Lemma \ref{claim3}, we   have $D(H) |X_i| = |V(H)|$, for every $i= 1, 2, \ldots , D(H)$, and so $D(H) ~|~|V(H)|$.\qed

\begin{theorem}\label{bestupperboung}
Let $G$ be a disconnected $d$-distinguishing  critical graph ($d\geq 5$) with $c(G)$ isomorphic connected components $H$ such that $c(G) < \frac{d}{2}$. If $\overline{G}\neq K_{d-1,d-1}, K_d$, then $c(G) \leq \frac{d-1}{3}$,  and so $D(H, D(H)) \leq  \frac{D(H)-3}{3}$. Also, $d\geq 7$ and   $D(H) \geq 6$.
\end{theorem}
\proof Since $\overline{G}$ is a connected  $d$-distinguishing  critical graph ($d\geq 5$), so $\overline{G}$ has a triangle as a induced subgraph, by Theorem \ref{triclawcrit}. Thus $G$ has a claw graph as a induced subgraph. Since $G$ is disjoint union of connected graphs $H$, so each of connected components $H$ has a claw. Hence, $\alpha (H) \geq 3$, and using $\alpha (G) = \alpha (H) c(G)$ and $\alpha (G) \leq d-1$  by Theorem \ref{characdiscritgra}, we get $3 c(G) \leq d-1$, and so $D(H, D(H)) \leq  \frac{D(H)-3}{3}$, because of $c(G) = D(H, D(H)) +1$ and $d = D(H)+1$ by Lemma \ref{lemineqcritic}. Finally, by  $3 c(G) \leq d-1$,  $c(G) \geq 2$ and $d = D(H)+1$, we get  $d\geq 7$ and $D(H) \geq 6$. \qed

Here, we can determine exactly all disconnected  $d$-distinguishing critical graphs for $d=5,6$.
\begin{corollary}~
\begin{itemize}
\item[{\rm (i)}] The only disconnected  $5$-distinguishing critical graphs are  $\overline{K_5}$ and $\overline{K_{4,4}}$.
\item[{\rm (ii)}] The only disconnected $6$-distinguishing critical graphs are $\overline{K_6}$  and $\overline{K_{5,5}}$.
\end{itemize}
\end{corollary}
\proof  It follows directly from Theorem \ref{bestupperboung}.\qed

Now we close this section by determining  the distinguishing critical trees. To do this, we need the following theorem:
\begin{theorem}{\rm \cite{EJC}}\label{uppboundisnumbtree}
If $T$ is a tree of order $n \geq 3$, then $D(T ) \leq \Delta (T)$. Furthermore, equality is achieved if
and only if $T$ is a symmetric tree or a path of odd length.
 \end{theorem}
\begin{theorem}
The only distinguishing critical tree is $K_2$.
\end{theorem}
\proof Let $T$ be a $d$-distinguishing critical tree of order $n\geq 3$. By Theorem \ref{uppboundisnumbtree},  we have $d\leq \Delta (T)$. Let $v$ be the vertex of $T$ with maximum degree $\Delta (T)$. Since $\Delta \geq d$ and $T$ does not have  any cycle, so we can make the star graph $K_{1,d}$ as a proper induced subgraph of $G$. Since $D(K_{1,d})=d$ and $T\neq K_{1,d}$ ($T$ is  a $d$-distinguishing critical tree), so we have a contradiction.\qed

\section{Conclusion}
In this paper we characterized the  $d$-distinguishing critical graphs. In particular,  we could determine exactly all  $d$-distinguishing critical graphs for $d=1,2, 3$. It can be seen that  all  such $d$-distinguishing critical graphs    are $k$-regular graph for some $k\leq d$. Also, we showed that the   disconnected    $d$-distinguishing critical graph with $c$ connected components such that $c\geq \frac{d}{2}$, is a regular graph. We propose the following conjecture:  

\begin{conjecture}~
\begin{itemize}
\item[(i)] If $G$ is a   $d$-distinguishing critical graph, then $G$ is a $k$-regular graph  for some $k\leq d$.
\item[(ii)] If $G$ is a  disconnected $d$-distinguishing critical graph, then each components of $G$ is a complete graph.
\end{itemize}
\end{conjecture}

We can extend the  definition of $d$-distinguishing critical graph to a  strong version.  A graph $G$ is called {\it strong $d$-distinguishing critical}, if $D(G)=d$ and $D(G-v)\neq D(G)$ for every vertex $v$ of $G$, where $G - v$ denotes the graph obtained from $G$ by removal of a vertex $v$ and
all edges incident to $v$. It is immediate that any   $d$-distinguishing critical graph is a  strong $d$-distinguishing critical graph, but  the converse is not true; the minimal asymmetric graphs are strong 1-distinguishing critical graphs which are not 1-distinguishing critical graphs.  This leads us to ask the following.
\begin{problem}
Characterize the  strong $d$-distinguishing critical graphs.
\end{problem}

\end{document}